\documentclass[12pt]{amsart}

\usepackage[margin = 1in]{geometry}
\usepackage{amsfonts, amsmath, amscd, amssymb,
latexsym, mathrsfs, mathtools, 
slashed, stmaryrd, verbatim, wasysym, bbm, url}
\usepackage[utf8]{inputenc}
\usepackage[backend=bibtex,bibstyle=authoryear,citestyle=alphabetic,giveninits=true,doi=false,isbn=false,url=false,backref]{biblatex}
\usepackage[usenames]{xcolor}
\usepackage{colortbl}
\usepackage{enumerate}
\usepackage[shortlabels]{enumitem}
\usepackage{lineno}
\linenumbers
\usepackage{scalerel}
\usepackage{soul}
\usepackage{diagbox}
\usepackage{array}
\usepackage{nicematrix}
\nolinenumbers
      
\DeclareFieldFormat{labelalphawidth}{\mkbibbrackets{#1}}

\defbibenvironment{bibliography}
  {\list
     {\printtext[labelalphawidth]{%
        \printfield{labelprefix}%
    \printfield{labelalpha}%
        \printfield{extraalpha}}}
     {\setlength{\labelwidth}{\labelalphawidth}%
      \setlength{\leftmargin}{\labelwidth}%
      \setlength{\labelsep}{\biblabelsep}%
      \addtolength{\leftmargin}{\labelsep}%
      \setlength{\itemsep}{\bibitemsep}%
      \setlength{\parsep}{\bibparsep}}%
      }
  {\endlist}
  {\item}
\DeclareNameAlias{author}{last-first}

\newtheorem*{theorem*}{Theorem}
\newtheorem{theorem}{Theorem}
\newtheorem{corollary}{Corollary}

\newtheorem{lemma}[corollary]{Lemma}
\newtheorem{proposition}[corollary]{Proposition}
\theoremstyle{definition}

\newtheorem{example}{Example}
\newtheorem{remark}[example]{Remark}

\numberwithin{equation}{section}
\let\oldsqrt\sqrt
\def\sqrt{\mathpalette\DHLhksqrt}
\def\DHLhksqrt#1#2{%
\setbox0=\hbox{$#1\oldsqrt{#2\,}$}\dimen0=\ht0
\advance\dimen0-0.2\ht0
\setbox2=\hbox{\vrule height\ht0 depth -\dimen0}%
{\box0\lower0.4pt\box2}}

\usepackage{verbatim} %% code from mathabx.sty and mathabx.dcl
\DeclareFontFamily{U}{mathx}{\hyphenchar\font45}
\DeclareFontShape{U}{mathx}{m}{n}{
      <5> <6> <7> <8> <9> <10>
      <10.95> <12> <14.4> <17.28> <20.74> <24.88>
      mathx10
      }{}
\DeclareSymbolFont{mathx}{U}{mathx}{m}{n}
\DeclareFontSubstitution{U}{mathx}{m}{n}
\DeclareMathAccent{\widecheck}{0}{mathx}{"71}

\renewcommand{\tilde}[1]{\widetilde{#1}}
\renewcommand{\bar}{\overline}

\renewcommand{\hat}[1]{\widehat{#1}}

\newcommand\eps\varepsilon
\renewcommand\epsilon\varepsilon

\newcommand{\redd}[1]{}

\newcommand{\abs}[1]{\left\lvert #1 \right\rvert}
\newcommand{\smallabs}[1]{\lvert #1 \rvert}

\newcommand{\norm}[1]{\lVert #1 \rVert}

\newcommand\supp{\operatorname{supp}}

\newcommand\paperintro%
        {%
         }
\newcommand\paperbody%
        {%
         }

\newcommand\bbN{\mathbb{N}}

\newcommand\bbR{\mathbb{R}}
\newcommand\bbS{\mathbb{S}}

\newcommand\cO{\mathcal{O}}

\DeclareMathAlphabet{\mathpzc}{OT1}{pzc}{m}{it}

\newcommand{\sbs}{\subset}

\usepackage[refpage]{nomencl}

\makenomenclature

\usepackage{tikz, tikz-cd, tkz-euclide}
\usetikzlibrary{calc,positioning,intersections,quotes,decorations.markings}
\makeatletter
\def\@tocline#1#2#3#4#5#6#7{\relax
  \ifnum #1>\c@tocdepth % then omitth
  \else
    \par \addpenalty\@secpenalty\addvspace{#2}%
    \begingroup \hyphenpenalty\@M
    \@ifempty{#4}{%
      \@tempdima\csname r@tocindent\number#1\endcsname\relax
    }{%
      \@tempdima#4\relax
    }%
    \parindent\z@ \leftskip#3\relax \advance\leftskip\@tempdima\relax
    \rightskip\@pnumwidth plus4em \parfillskip-\@pnumwidth
    #5\leavevmode\hskip-\@tempdima
      \ifcase #1
       \or\or \hskip 1em \or \hskip 2em \else \hskip 3em \fi%
      #6\nobreak\relax
    \hfill\hbox to\@pnumwidth{\@tocpagenum{#7}}\par% <---- \dotfill -> \hfill
    \nobreak
    \endgroup
  \fi}
\makeatother

\def\annu#1{_{% 
  \vbox{\hrule height .2pt 
    \kern 1pt 
    \hbox{$\scriptstyle {#1}\kern 1pt$}% 
  }\kern-.05pt 
  \vrule width .2pt 
}}

\allowdisplaybreaks
\usepackage[colorlinks,citecolor=blue,urlcolor=red,bookmarks=false,hypertexnames=true]{hyperref} 

\makeatletter
\def\keywords{\xdef\@thefnmark{}\@footnotetext}
\makeatother

\title{A Rellich-type inequality for the sphere}
\author{Ali Barki}
\address{UPL, Univ. Paris Nanterre, CNRS, F92000 Nanterre France, and Department of Mathematics, Imperial College London, 180 Queen's Gate, SW7 2AZ, London, United Kingdom}
\email{ali.barki@ens-rennes.fr}
\thanks{The first author was supported by the Labex MME-DII, funded by ANR-11-LBX-0023-01.}

\author{Shubham Gupta}
\address{Department of Mathematics, Technion--Israel Institute of Technology, Haifa 32000 Israel}
\email{sg1019@campus.technion.ac.il}
\thanks{The second author was supported by the ISF-DFG Lead Agency research
grant (ISF No. 1531/25).}

\author{Yaozhong Qiu}
\address{UPL, Univ. Paris Nanterre, CNRS, F92000 Nanterre France}
\email{yqiu@parisnanterre.fr}
\thanks{The third author was supported by the President's PhD Scholarship of Imperial College London. This project has received funding from the European Union’s Horizon 2020 research and innovation programme under the Marie Sk\l{}odowska-Curie grant agreement No 101034255. \scalerel*{\includegraphics{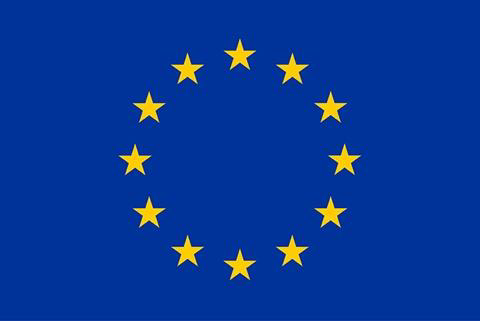}}{A}}

\begin{document}
\begin{abstract}
We prove a Rellich-type inequality of the form 
$$\int \abs{(\Delta + 1)u}^2dx \geq C_d \int \frac{\abs{u}^2}{(1 - \abs{x}^2)^2}dx $$
for smooth functions $u$ compactly supported in the unit ball in $\mathbb{R}^d$ for $d \geq 1$. The obtained constant is asymptotically sharp up to the first three orders, as $d \rightarrow \infty$. We also show that in the complement of the ball, the inequality holds only in dimension one, and does not hold true with a positive $C_d>0$ for $d \geq 2$. The motivation is that if $\hat{u}$ denotes the Fourier transform of $u$, then the inequality is associated with the uncertainty principle
$$ \left(\int \abs{u}^2dx\right)^2 \lesssim \int (1-\abs{x}^2)^2\abs{u}^2dx \int (1 - \abs{\xi}^2)^2\smallabs{\hat{u}}^2d\xi,$$
which says that a function cannot be simultaneously concentrated on a sphere in both space and frequency. This is in analogy with the classical Hardy inequality together with the Heisenberg uncertainty principle. 
\end{abstract}

\keywords{2020 \emph{Mathematics Subject Classification.} Primary 26D10, 35A23}
\keywords{\emph{Keywords.}
Hardy inequality, Rellich inequality, uncertainty principles, spherical singularity.}

\maketitle

\section{Introduction and main results}
We begin by recalling the classical Hardy inequality 
\begin{equation}\label{eq:classical-hardy} 
    \int \abs{\nabla u}^2dx \geq \left(\frac{d-2}{2}\right)^2 \int \frac{\abs{u}^2}{\abs{x}^2}dx
\end{equation}
valid for functions $u$ compactly supported outside the origin. It was first established in the one-dimensional setting by Hardy \cite{hardy1920note} and later extended to higher dimensions \cite{leray1933etude}. By the Cauchy-Schwarz inequality, it implies the well-known  Heisenberg uncertainty principle 
\begin{align*}
    \left(\frac{d-2}{2}\right)^2 \left(\int \abs{u}^2dx\right)^2 &\leq \int \abs{x}^2\abs{u}^2dx \int \abs{\nabla u}^2dx \\
    &= \int \abs{x}^2\abs{u}^2dx \int \abs{\xi}^2\smallabs{\hat{u}}^2d\xi,
\end{align*}
where $\hat{u}$ denotes the Fourier transform of $u$.
Although the inequality is not sharp, owing to the nonoptimal use of Hardy inequality, the underlying uncertainty principle holds with the optimal constant $(d/2)^2\geq ((d-2)/2)^2$. Nevertheless, the inequality captures an important uncertainty principle: a function $u$ cannot be simultaneously concentrated at two points $x \in A = \{0\}$ and $\xi \in B = \{0\}$ in space and frequency respectively. Moreover, by translation and modulation, it can be shown $A$ and $B$ can be any pair of singletons, not necessarily zero.   

The question we would like to address is whether similar statements can be made for a more general class of sets $A, B \sbs \bbR^d$, which are not necessarily singletons. In this paper, we specifically consider the case when both $A$ and $B$ are spheres, which leads us to the following question:

\textbf{Question.}
    Let $d\in \mathbb{N}$ and $\mathbb{S}^{d-1}$ be the unit sphere in $\mathbb{R}^d$. Does there exist a positive constant $C_d>0$  such that the following Rellich-type inequality holds 
    \begin{equation}\label{eq:spherical-rellich}
    \int \abs{(\Delta + 1)u}^2dx \geq C_d \int \frac{\abs{u}^2}{(1 - \abs{x}^2)^2}dx,
    \end{equation}
    for every smooth function $u$ compactly supported in $\mathbb{R}^d \setminus \mathbb{S}^{d-1}$?

In the same way the Hardy inequality \eqref{eq:classical-hardy} provides a route to the uncertainty principle when $A$ and $B$ are singletons, \eqref{eq:spherical-rellich} naturally leads to uncertainty principle for spheres. Indeed, since 
$$
\widehat{(\Delta+1)u}(\xi) = (1-|\xi|^2) \widehat{u}(\xi),
$$
the validity of \eqref{eq:spherical-rellich} would imply, via Cauchy--Schwarz and Plancherel's theorem
\begin{equation}\label{eq:spherical-uncertainty}
    C_d \left(\int \abs{u}^2dx\right)^2 \leq \int (1 - \abs{x}^2)^2\abs{u}^2dx \int (1 - \abs{\xi}^2)^2\smallabs{\hat{u}}^2d\xi. 
\end{equation}
The validity of such an inequality would therefore assert that a function $u$ cannot be simultaneously concentrated near the unit sphere in both space and frequency. 

The goal of this paper is therefore to provide a preliminary investigation into the validity of such uncertainty principles from the perspective of Hardy and Rellich type inequalities and to provide an answer in the spherical case.  Since the sphere separates the space $\mathbb{R}^d$ into two regions, it is natural to examine \eqref{eq:spherical-rellich} on each of them. As we shall see, this turns out to be decisive: the inequality holds true inside the unit ball in all dimensions, while outside, it holds only in dimension one.

We now state our main results. Let $B$ denote the open unit ball centered at the origin $$B = \{x \in \bbR^d : \abs{x} < 1\},$$ 
and let $B^c$ denote the interior of its complement. 

\begin{theorem}\label{thm:int}
Let $d \geq 1$ and $u \in C_c^\infty(B)$. Then we have
    \[ \int_B \abs{(\Delta + 1)u}^2dx \geq C_d \int_B \frac{\abs{u}^2}{(1 - \abs{x}^2)^2}dx,\]
    where
    $$C_1 = 3, \quad and \quad C_d = 8d^2 + 16 d - 32 \; for \; d \geq 2.$$ 
    Moreover the constant is asymptotically sharp, that is, 
    $$
    C_d^\textnormal{opt} = 8d^2 + 16 d - 32 + O(d^{-1}), \quad \operatorname{as} \; d \to \infty,
    $$
    where $C_d^\textnormal{opt}$ denotes the optimal constant in the inequality.
\end{theorem}

\begin{theorem}\label{thm:ext}
Let $u \in C_c^\infty(B^c)$. Then in dimension $d=1$, we have 
    \[ \int_{B^c} \abs{(\Delta + 1)u}^2dx \geq \tilde{C}_1 \int_{B^c} \frac{\abs{u}^2}{(1 - \abs{x}^2)^2}dx \]
    with $\tilde{C}_1 = 3$, and otherwise for $d \geq 2$, such an inequality does not hold for any $\tilde{C}_d > 0$.  
\end{theorem}

The failure of the Rellich inequality for $d \geq 2$ is established using test functions involving spherical harmonics of arbitrary high degree. The obstruction disappears when one restricts to radial functions, since they only involve degree zero harmonics. Indeed, as we shall show in Proposition \ref{prop:radial-functions}, one can recover the Rellich inequality in the whole space $\mathbb{R}^d \setminus \mathbb{S}^{d-1}$ for radial functions. Taken together, Theorems \ref{thm:int} and \ref{thm:ext} answer the question posed at the beginning: the Rellich type inequality \eqref{eq:spherical-rellich} holds true for all functions $u \in C_c^\infty(\mathbb{R}^d \setminus \mathbb{S}^{d-1})$ if and only if $d=1$. Nevertheless, Theorem \ref{thm:int} and Proposition \ref{prop:radial-functions} establishes the spherical uncertainty principle for functions supported inside the unit ball and for radial functions, respectively. 

Before concluding this section, let us place our results in a broader context. At least from a philosophical point of view, the question of understanding uncertainty principles for general sets $A,B\subset\mathbb{R}^d$ is not new and has been addressed in several works, see for instance \cite{donoho1989uncertainty, jaming2007nazarov, wigderson2021uncertainty}. However, the quantitative nature of the results in these works is different from that of the results established here. Our approach is based on Hardy and Rellich type inequalities, a perspective that, to the best of our knowledge, has not been explored before in this context.

Concerning the specific inequality associated with the sphere studied here, the weight $(1-|x|^2)^{-2}$ is closely related to distance from the boundary of the domain. Rellich inequalities for such weights have been studied extensively in the literature, see for instance \cite{owen1999hardy, barbatis2006class}. The novelty of our results lies in the use of a different operator $\Delta+1$, which is precisely the operator needed to establish the spherical uncertainty principle. From this perspective, our work also raises a natural question concerning the class of operators that can be used in place of $\Delta+1$, leading to suitable uncertainty statements. As far as we are aware, such questions have not been studied in a systematic way. This work is only a first step towards that direction, and we hope that this general curiosity leads to some interesting mathematics.

For the interested reader, the literature on Hardy and Rellich type inequalities is rich and extremely vast and we cannot hope to provide an exhaustive list of the many developments in the study thereof. Nevertheless, we would like to direct the reader to the following literature which we believe are particularly relevant to our results and methods \cite{beckner1995pitt, yafaev1999sharp, gesztesy2018factorizations, cazacu2022sharp}.
\section{Proof of main results}
We start with a preliminary Hardy-type inequality which is needed to control a remainder term appearing in the Rellich expansion of squares. To ease notation, here and in the sequel, we set $r = \abs{x}$. We also assume $u$ is real-valued without loss of generality.

\begin{lemma}\label{lem:spherical-hardy}
    Assume $u \in C_c^\infty(B)$. Then it holds 
    \[ 2d(n+1) \int \frac{\abs{u}^2}{(1 - r^2)^{n+1}}dx \leq \int \frac{\abs{\nabla u}^2}{(1 - r^2)^n}dx \]
    for all $n \geq 0$. 
\end{lemma}

\begin{proof}
    We expand the square 
    \[ \int \abs{(1-r^2)^{-n/2}\nabla u + C(1-r^2)^{-n/2-1}xu}^2dx \]
    where $x$ is regarded as the vector $(x_1, \cdots, x_d)$. The square terms are 
    \[ \int \frac{\abs{\nabla u}^2}{(1 - r^2)^n}dx, \quad C^2 \int \frac{\abs{x}^2\abs{u}^2}{(1 - r^2)^{n+2}}dx \]
    while the cross term can be studied by integrating by parts: 
    \begin{align*}
        2C \int \frac{u(x \cdot \nabla u)}{(1 - r^2)^{n+1}}dx &= 2C \sum_{j=1}^d \frac{u x_ju_j}{(1-r^2)^{n+1}}dx \\
        &= C\sum_{j=1}^d \int \frac{x_j(u^2)_j}{(1-r^2)^{n+1}}dx \\
        &= -C\sum_{j=1}^d \int \frac{u^2}{(1-r^2)^{n+1}}dx + C(n+1)\sum_{j=1}^d \int \frac{-2x_j^2u^2}{(1-r^2)^{n+2}}dx \\
        &= -Cd \int \frac{u^2}{(1-r^2)^{n+1}}dx - 2C(n+1) \int \frac{\abs{x}^2u^2}{(1-r^2)^{n+2}}dx.
    \end{align*}
    By taking $C = 2(n+1)$ we obtain the claimed result by rearranging 
    \[ \int \frac{\abs{\nabla u}^2}{(1 - r^2)^n}dx + 2C \int \frac{u(x \cdot \nabla u)}{(1 - r^2)^{n+1}}dx + C^2 \int \frac{\abs{x}^2\abs{u}^2}{(1 - r^2)^{n+2}}dx \geq 0. \qedhere \]
\end{proof}

Note the lemma describes the uncertainty principle that a function $u$ compactly supported in $B$ cannot be simultaneously concentrated on a sphere $A = \bbS^{d-1}$ and a point $B = \{0\}$ in space and frequency respectively. 

To prove Theorem \ref{thm:int}, we provide a proof of the following expansion (valid in both the interior and exterior of the ball) of the square.

\begin{lemma}\label{lem:rellich-expansion}
    Let $f, g$ be two smooth radial functions and assume $u$ is smooth and compactly supported in either $B$ or $B^c$. Then it holds 
    \begin{equation}\label{eq:rellich-expansion}
    \begin{aligned}
        \int \abs{(\Delta + 1)u}^2dx &\geq \int \left(2\frac{f'}{r} - f^2\right)\abs{x \cdot \nabla u}^2dx \\
        &+ \int \left((2-d)f - rf' + 2g\right)\abs{\nabla u}^2dx \vphantom{\int \abs{x \cdot \nabla u}^2\left(2\frac{f'}{r} - f^2\right)dx} \\
        &+ \int \left(-g^2 + rf' + df - \Delta g - 2g + r(fg)' + dfg\right)\abs{u}^2dx. \vphantom{\int \abs{x \cdot \nabla u}^2\left(2\frac{f'}{r} - f^2\right)dx}
    \end{aligned}
    \end{equation}
\end{lemma}

\begin{proof}
    We expand the square 
    \[ \int \abs{(\Delta + 1)u + f(x \cdot \nabla u) + gu}^2dx \geq 0. \]
    The three squares are 
    \[ \int \abs{(\Delta + 1)u}^2dx, \quad \int f^2\abs{x \cdot \nabla u}^2dx, \quad \int g^2u^2dx. \] 
    For the first cross term, we have two components. The first contribution is 
    \begin{align*}
         2\int f(\Delta u)(x \cdot \nabla u)dx &= -2\int (\nabla f \cdot \nabla u)(x \cdot \nabla u)dx - 2\int f\nabla(x \cdot \nabla u) \cdot \nabla udx 
         \vphantom{-2\int \frac{f'}{r}\abs{x \cdot \nabla u}^2dx - 2\sum_{j=1}^d \int fu_j^2dx + \sum_{i, j = 1}^d \int \left(\frac{f'}{r}x_i^2 + f\right)u_j^2dx} \\
         &= -2 \int \frac{f'}{r} \abs{x \cdot \nabla u}^2dx - 2 \sum_{i, j = 1}^d \int f\partial_j(x_iu_i)u_jdx 
         \vphantom{-2\int \frac{f'}{r}\abs{x \cdot \nabla u}^2dx - 2\sum_{j=1}^d \int fu_j^2dx + \sum_{i, j = 1}^d \int \left(\frac{f'}{r}x_i^2 + f\right)u_j^2dx} \\ 
         &= -2\int \frac{f'}{r}\abs{x \cdot \nabla u}^2dx - 2\sum_{i, j = 1}^d \int f(x_iu_{ij} + u_i\delta_{ij})u_jdx 
         \vphantom{-2\int \frac{f'}{r}\abs{x \cdot \nabla u}^2dx - 2\sum_{j=1}^d \int fu_j^2dx + \sum_{i, j = 1}^d \int \left(\frac{f'}{r}x_i^2 + f\right)u_j^2dx} \\ 
         &= -2\int \frac{f'}{r}\abs{x \cdot \nabla u}^2dx - 2\sum_{j=1}^d \int fu_j^2dx - \sum_{i, j = 1}^d \int fx_i(u_j^2)_idx 
         \vphantom{-2\int \frac{f'}{r}\abs{x \cdot \nabla u}^2dx - 2\sum_{j=1}^d \int fu_j^2dx + \sum_{i, j = 1}^d \int \left(\frac{f'}{r}x_i^2 + f\right)u_j^2dx} \\ 
         &= -2\int \frac{f'}{r}\abs{x \cdot \nabla u}^2dx - 2\sum_{j=1}^d \int fu_j^2dx + \sum_{i, j = 1}^d \int \left(\frac{f'}{r}x_i^2 + f\right)u_j^2dx.  
    \end{align*}
    The final addend is 
    \[ \sum_{i, j = 1}^d \int \left(\frac{f'}{r}x_i^2 + f\right)u_j^2dx = \int rf'\abs{\nabla u}^2dx + d\int f\abs{\nabla u}^2dx \]
    so that we find 
    \begin{equation}\label{eq:expansion-1.1} 
        2\int  f(\Delta u)(x \cdot \nabla u)dx = \left((d - 2)f + rf'\right)\int \abs{\nabla u}^2dx - 2\int \frac{f'}{r}\abs{x \cdot \nabla u}^2dx. 
    \end{equation}
    The second contribution is 
    \begin{align}
         2\int fu(x \cdot \nabla u)dx &= 2 \sum_{i=1}^d \int fx_iu_iudx 
         \vphantom{-\sum_{i=1}^d \int \frac{f'}{r}x_i^2u^2dx - \sum_{i=1}^d \int fu^2dx} \nonumber \\ 
         &= \sum_{i=1}^d \int fx_i(u^2)_idx 
         \vphantom{-\sum_{i=1}^d \int \frac{f'}{r}x_i^2u^2dx - \sum_{i=1}^d \int fu^2dx} \nonumber \\
         &= -\sum_{i=1}^d \int \frac{f'}{r}x_i^2u^2dx - \sum_{i=1}^d \int fu^2dx \nonumber \\ 
         &= -\int rf'u^2dx - d\int fu^2dx 
         \vphantom{-\sum_{i=1}^d \int \frac{f'}{r}x_i^2u^2dx - \sum_{i=1}^d \int fu^2dx} \label{eq:expansion-1.2}
    \end{align}
    and together this completes the expansion of the cross term $\int 2f((\Delta + 1)u)(x \cdot \nabla u)dx$. 

    The second cross term is 
    \begin{align}
         2\int f g u(x \cdot \nabla u)dx &= \int f g x \cdot (\nabla u^2)dx
         \vphantom{-\sum_{i=1}^d \int \frac{f'}{r}x_i^2u^2dx - \sum_{i=1}^d \int fu^2dx} \nonumber 
         \vphantom{-\sum_{i=1}^d \int \frac{f'}{r}x_i^2u^2dx - \sum_{i=1}^d \int f g u^2dx} \nonumber \\ 
         &= \sum_{i=1}^d \int f g x_i(u^2)_idx 
         \vphantom{-\sum_{i=1}^d \int \frac{(fg)'}{r}x_i^2u^2dx - \sum_{i=1}^d \int f g u^2dx} \nonumber \\
         &= -\sum_{i=1}^d \int \frac{(fg)'}{r}x_i^2u^2dx - \sum_{i=1}^d \int f gu^2dx \nonumber \\ 
         &= -\int r(fg)'u^2dx - d\int fg u^2dx
         \vphantom{-\sum_{i=1}^d \int \frac{f'}{r}x_i^2u^2dx - \sum_{i=1}^d \int fu^2dx}. \label{eq:expansion-2}
    \end{align}
    while the third cross term is 
    \begin{align}
        2\int ((\Delta + 1)u)(gu)dx &= \int g\Delta(u^2)dx - 2\int g\abs{\nabla u}^2dx + 2\int gu^2dx \nonumber \\
        &= \int \Delta(g)u^2dx - 2\int g\abs{\nabla u}^2dx + 2\int gu^2dx. \label{eq:expansion-3} 
    \end{align}
    Rearranging the identities \eqref{eq:expansion-1.1}, \eqref{eq:expansion-1.2}, \eqref{eq:expansion-2}, \eqref{eq:expansion-3} yields \eqref{eq:rellich-expansion}. 
\end{proof}

\begin{remark}
    One may also expand instead the square 
    \[ \int \abs{(\Delta + 1)u + \nabla f_0 \cdot \nabla u + gu}^2dx \geq 0 \]
    for some other choice $f_0$; but since $\nabla f_0(x) = \frac{x}{r}f_0'(r)$, we can write $\nabla f_0 = fx$ for $f = f_0'/r$. 
\end{remark}

Following the proof of the Rellich inequality by expansion of the squares, see for instance \cite{gesztesy2018factorizations}, we now choose a pair of suitable $f, g$ such that the right hand side of \eqref{eq:rellich-expansion} is reminiscent of the expected right side in the statement of the two theorems. Inspired by the Rellich setting, we take 
\[ f(r) = \frac{a}{1-r^2}, \quad g = \frac{b}{1-r^2} + \frac{c}{(1 - r^2)^2} \]
for $a, b, c \in \bbR$ to be determined. 

The second term in $g$ may appear artificial but is motivated by the observation $2f'(r)/r - f(r)^2 = a(4-a)(1 - r^2)^{-2}$ so that if $a \in \bbR \setminus (0, 4)$ and $a(4-a) \leq 0$ then Cauchy-Schwarz implies 
\[ \int \left(2\frac{f'}{r} - f^2\right)\abs{x \cdot \nabla u}^2dx \geq a(4-a) \int \frac{r^2}{(1-r^2)^2}\abs{\nabla u}^2dx \]
and thus can be absorbed into the weighted gradient term in \eqref{eq:rellich-expansion}. Then 
\[ (2-d)f(r) - rf'(r) + \frac{2b}{1-r^2} + r^2\left(2\frac{f'(r)}{r} - f(r)^2\right) = \frac{(2a+2b-da) - (a^2+2b-da)r^2}{(1 - r^2)^2} \]
which is proportional to $(1 - r^2)^{-1}$ if $2a+2b-da = a^2+2b-da$. This occurs when $a \in \{0, 2\}$ and hence $a = 0$ since we assumed $a \in \bbR \setminus (0, 4)$. This yields 
\[ \int \abs{(\Delta + 1)u}^2dx \geq 2b \int \frac{\abs{\nabla u}^2}{1 - r^2}dx - b \int \frac{(2 + b + 2d - (-4 + b + 2d)r^2 + 2r^4)}{(1 - r^2)^3}u^2dx. \]
But as $r \rightarrow 1$, the numerator on the weighted $L^2$-term approaches $-8b$ and we cannot make both weights simultaneously nonnegative near $r = 1$. On the other hand, if $a \in (0, 4)$ and $a(4-a) > 0$ then we may bound the weighted radial derivative term in \eqref{eq:rellich-expansion} by zero. Then 
\[ (2 - d)f(r) - r'f(r) + \frac{2 b}{1 - r^2} = \frac{(2a + 2b - da) - (4a + 2b - da)r^2}{(1 - r^2)^2} \]
which implies again $a = 0$. The second term in $g$ therefore appears in order to remedy this problem by allowing for some additional cancellation at the cost of complicating the weight on the $L^2$-term.

\begin{proof}[Proof of Theorem \ref{thm:int} for $d \geq 2$]
    Assuming once again $c \in \bbR$ and $a \in \bbR \setminus (0, 4)$ we have 
    \begin{align*}
        \omega_1(r) &= (2 - d)f(r) - r f'(r) + 2g(r) + r^2 \left(2 \frac{f'(r)}{r} - f(r)^2\right) \\
        &= \frac{(2(a + b + c) - da) - (a^2 + 2b - da)r^2}{(1 - r^2)^2} \vphantom{(2 - d)f(r) - r f'(r) + 2g(r) + r^2 \left(2 \frac{f'(r)}{r} - f(r)^2\right)} 
    \end{align*}
    and to make the numerator proportional to $1 - r^2$ it suffices to send $r \rightarrow 1$ in the numerator and solve for zero (which is equivalent to equating the two coefficients); this implies $c = \frac{1}{2}(a^2 - 2a)$ which yields 
    \begin{equation}\label{eq:omega-1}
        \omega_1(r) = \frac{a^2+2b-da}{1-r^2}.
    \end{equation}
    The weight on the $L^2$-term, which we denote by $\omega_2(r)$, is then 
    \begin{equation}\label{eq:omega-2}
        \omega_2(r) = \frac{C_1 + C_2r^2 + C_3r^4 + C_4r^6}{4(1 - r^2)^4}
    \end{equation}
    for some $C_1, C_2, C_3, C_4 \in \bbR$ depending on $a, b, d$, specifically
    \begin{align*}
        C_1 &= -a^4 + 2a^3(2 + d) - 4b(2 + b + 2d) - 4a^2(2 + b + 3d) + 4a(2 + 5d + b(2 + d)) \\ 
        C_2 &= a(88 - 8b(-1 + d)- 28d) - 2a^3(-6 + d) + 8b(-1 + b + 2d) + 4a^2(-16 + b + 3d) \\
        C_3 &= -4(a^2 + a(2 - b(-4 + d) - 3d) + b(-2 + b + 2d)) \\
        C_4 &= 8(a + b) - 4da. 
    \end{align*}
    As before, we send $r \rightarrow 1$ in the numerator and solve $C_1 + C_2 + C_3 + C_4 = 0$ for $a$. That this factors $1-r^2$ and not simply $1-r$ is since only even powers of $r$ appear in the numerator (and thus the numerator is even) so the second root of $1-r^2$ at $r = -1$ is automatically factored by the same condition. We deduce $a \in (\bbR \setminus (0, 4)) \cap \{0, 2, 6, 8\} = \{0, 6, 8\}$. The case $a = 0$ implies $c = 0$ and can be ruled out by the previous discussion. 
    
    The remaining two cases $a \in \{6, 8\}$ are analysed separately. It turns out that the correct choice is $a = 8$, in which case \eqref{eq:omega-2} becomes 
    \[ \frac{(-b(50+b-6d)+104(-6+d)) + (64+28b+b^2-2(8+3b)d)r^2 + (-2(8+b-4d))r^4}{(1 - r^2)^3} \]
    and repeating the previous argument (sending $r \rightarrow 1$ in the numerator and solving for zero for $b$) yields $b = 4(-6+d)$ and thus together with \eqref{eq:omega-1} it holds
    \[ \int_B \abs{(\Delta + 1)u}^2dx \geq 16 \int_B \frac{\abs{\nabla u}^2}{1 - r^2}dx + 8 \int_B \frac{d(d-6) - 4}{(1 - r^2)^2}u^2dx \] 
    since $u \in C_c^\infty(B)$ and hence $r \leq 1$. By applying Lemma \ref{lem:spherical-hardy}, we obtain the desired inequality with the expected constant 
    \[ C_d = 16 \times 2d(1+1) + 8(d(d-6)-4) = 8(d^2 + 2d - 4). \]
    On the other hand, if $a = 6$, it can be verified the same arguments yield (the details of which are omitted for brevity) the same inequality with constant $3(d(d-4)-16)$ which is positive only if $d \geq 7$. If $a \in (0, 4)$ then $c = a$ and $a \in \{0, \frac{24}{5}\}$, which is inconclusive. 

    For the asymptotic optimality, a straightforward computation with the test function 
    \[ u_*(x) = u_*(r) = (1 - r^2)^2 \]
    leads to 
    \begin{align*}
        \dfrac{\int_B \abs{(\Delta + 1)u_*}^2dx}{\int_B \dfrac{u_*^2}{(1 - r^2)^2}dx} &= \frac{\int_0^1 r^{d-1}((1 - 4d) + (4d + 6)r^2 + r^2)^2dr}{\int_0^1 r^{d-1}(1 - r^2)^2dr } \\
        &= \frac{8(d^4 + 16d^3 + 72d^2 + 64d + 6)}{(d + 6)(d + 8)} \vphantom{\dfrac{\int_B \abs{(\Delta + 1)u_*}^2dx}{\int_B \dfrac{u_*^2}{(1 - r^2)^2}dx}}
    \end{align*}
    and it is readily checked this rational function of $d$ is equal to $C_d + O(1/d)$ for large $d$. Of course, $u_*$ does not belong to $C_c^\infty(B)$. Nevertheless, the integration by parts arguments used in the proof of Lemma \ref{lem:spherical-hardy} hold whenever the boundary terms vanish and for this it suffices that $u \in H_0^2(\Omega)$, and $u_* \in H_0^2(\Omega)$ since $u \in C^\infty(\bar{B})$ and $u_*$ and $\nabla u_*$ both vanish at the boundary $\partial B$. 
\end{proof}

Having completed the proof in dimension $d \geq 2$, we now proceed to the case of dimension $1$ where the weighted radial derivative $\abs{x \cdot \nabla u}^2$ is simply $r^2\abs{\nabla u}^2$, and which permits a simplification (in fact in all dimensions $d \leq 3$ assuming $u$ is radial). 

\begin{proof}[Proof of Theorems \ref{thm:int} and \ref{thm:ext} for $d = 1$]
    Assuming $u$ is radial, the inequality \eqref{eq:rellich-expansion} simplifies, for the same choice of $f$ and $g$, to the case $a \in \bbR \setminus (0, 4)$ in the previous proof, except we permit $a \in \bbR$. Thus we may take $a = 2$ in which case $c = \frac{1}{2}(-2a + a^2) = 0$. This yields a double root in the numerator of the analogous $\omega_2$ and we obtain 
    \[ \int \abs{(\Delta + 1)u}^2dx \geq 2(2+b-d) \int \frac{\abs{\nabla u}^2}{1-r^2}dx + \int \frac{(-b^2-2b+2d) + 2(2+b-d)r^2}{(1 - r^2)^2}u^2dx. \]
    The choice $b = d-2$ annihilates the first addend and the coefficient of $r^2$ in the numerator of the second and thus leads to the constant $d(4-d)$ which is equal to $3$ for $d = 1$. This also provides a partial answer in the case of radial functions compactly supported in $B^c$ for $d \in \{2, 3\}$. 
\end{proof} 

\begin{proof}[Proof of Theorem \ref{thm:ext} for $d \geq 2$]
We construct a sequence $(u_n)_{n\geq1}\subset C_c^\infty(B^c)$ such that
\begin{equation}\label{eq:zero-limit}
    \lim_{n \to \infty} \dfrac{\int_{B^c} \abs{(\Delta + 1)u_n}^2dx}{\int_{B^c} \dfrac{\abs{u_n}^2}{(r^2 - 1)^2}dx} = 0. 
\end{equation}
Let $n \in \bbN$ and $Y_n$ be a spherical harmonic of degree $n$, normalized in $L^2(\bbS^{d-1})$. Writing $u(x) = v(r) Y_n(\theta)$ in the polar coordinates and decomposing the Laplacian into radial and spherical components, we get
\[ (\Delta + 1)u(x) = v''(r)Y_n(\theta) + \frac{d-1}{r}v'(r)Y_n(\theta) + \frac{1}{r^2}v(r)\Delta_{\bbS^{d-1}}Y_n(\theta) + v(r)Y_n(\theta). \]
Since $Y_n$ is an eigenfunction of $-\Delta_{\bbS^{d-1}}$ with eigenvalue $n(n+d-2)$, the expression simplifies to 
\[ (\Delta + 1)u(x) = \left(v''(r) + \frac{d-1}{r}v'(r) - \frac{n(n+d-2)}{r^2}v(r) + v(r)\right) Y_n(\theta). \]
Making the change of variable $v = r^{-(d-1)/2}w$, we have
\begin{equation}\label{eq:delta+1-expression}
    (\Delta + 1)u(x) = r^{-\frac{d-1}{2}}\left( w''(r) + \left(1- \frac{\nu_n^2- \frac{1}{4}}{r^2}\right)w(r)\right)  Y_n(\theta),    
\end{equation}
where
\[ \nu_n := n + \frac{d-2}{2}. \]
The equation on the right hand side is a variation of the familiar Bessel's equation, a solution of which is $w_n = \sqrt{r}J_{\nu_n}(r)$, where $J_{\nu_n}$ is the first Bessel function of order $\nu_n$. Note that while this choice of $w$, and hence $u$, makes the numerator of \eqref{eq:zero-limit} zero, it is not compactly supported. Hence, appropriate cutoffs are needed to make these functions admissible. To this end, we define two cutoffs, $\chi$ and $\eta$ around $r = 1$ and $r = +\infty$ respectively by
\[ \chi(r) = \begin{cases} 0, &0 \leq r \leq 2, \\ 1, &r \geq 3, \end{cases}\]
and
\[ \eta(r) = \begin{cases} 1, &0 \leq r \leq 1, \\ 0, &r \geq 2. \end{cases} \]
Let $(R_n)_{n \geq 1}$ be a divergent sequence $R_n \rightarrow +\infty$ to be determined later. Let $a_n(r) = \chi(r)\eta(r/R_n)$ and define $u_n(r, \theta) = r^{-\frac{d-1}{2}} a_n(r) w_n (r) Y_n(\theta)$. Note that $u_n \in C_c^\infty(B^c)$ and  
\[ \supp(u_n) \sbs \{x \in \bbR^d \mid 2 \leq \abs{x} \leq 2R_n\}.\]
The expression \eqref{eq:delta+1-expression} along with the fact that $w_n$ solves
\[ w''(r) + \left(1 - \frac{\nu_n^2 - \frac{1}{4}}{r^2}\right)w(r) = 0 \]
implies 
\[ (\Delta + 1)u_n = r^{-\frac{d-1}{2}}\left(2a_n'(r)w_n'(r) + a_n''(r)w(r)\right)Y_n(\theta). \]
We start by estimating the numerator of \eqref{eq:zero-limit}. Expressing in polar coordinates and ignoring the surface measure of the sphere, we get 
\begin{equation}\label{eq:delta+1-un}
    \int_{B^c} \abs{(\Delta + 1)u_n}^2dx \propto \int_1^\infty \abs{2a_n'(r)w_n'(r) + a_n''(r)w_n(r)}^2dr.
\end{equation}
By differentiating, we obtain 
\begin{align*}
    a_n'(r) &= \chi'(r)\eta\left(\frac{r}{R_n}\right) + \frac{\chi(r)}{R_n}\eta'\left(\frac{r}{R_n}\right), \\ 
    a_n''(r) &= \chi''(r)\eta\left(\frac{r}{R_n}\right) + 2\frac{\chi'(r)}{R_n}\eta'\left(\frac{r}{R_n}\right) + \frac{\chi(r)}{R_n^2}\eta''\left(\frac{r}{R_n}\right)
\end{align*}
and we note it follows from the properties of the cutoffs that 
\[ \supp(a_n') \cup \supp(a_n'') \sbs [2, 3] \cup [R_n, 2R_n]. \]

In order to understand the decay of $w_n$ and $w_n'$ on $[2, 3]$, we recall the following power series representation 
\begin{equation}\label{eq:bessel-representation}
    J_\nu(r) = \sum_{k = 0}^\infty \frac{(-1)^k}{k!\Gamma(\nu + k + 1)}\left(\frac{r}{2}\right)^{2k+\nu}
\end{equation}
for the Bessel function, see for instance \cite{abramowitz1948handbook, NIST:DLMF}. For $\nu \geq 1$ and $k \in \bbN$, since $\Gamma(\nu + k + 1) \geq \Gamma(\nu + 1)(\nu + 1)^k$, we have the estimate
\begin{equation}\label{eq:J-upper-bound}
\begin{aligned}
    \sup_{2 \leq r \leq 3} \abs{J_\nu(r)} &\leq \sup_{2 \leq r \leq 3} \frac{(r/2)^\nu}{\Gamma(\nu+1)}\sum_{k=0}^\infty \frac{1}{k!}\left(\frac{r^2}{4\Gamma(\nu+1)}\right)^k \\
    &= \sup_{2 \leq r \leq 3} \frac{(r/2)^\nu}{\Gamma(\nu+1)}\exp\left(\frac{r^2}{4(\nu+1)}\right) \vphantom{\frac{(r/2)^\nu}{\Gamma(\nu+1)}\sum_{k=0}^\infty \frac{1}{k!}\left(\frac{r^2}{4\Gamma(\nu+1)}\right)^k} \\
    &\leq e^{9/8} \frac{(3/2)^\nu}{\Gamma(\nu+1)} \vphantom{\frac{(r/2)^\nu}{\Gamma(\nu+1)}\sum_{k=0}^\infty \frac{1}{k!}\left(\frac{r^2}{4\Gamma(\nu+1)}\right)^k}.
\end{aligned}
\end{equation}
Together with the recurrence relation $J_\nu'(r) = \frac{\nu}{r}J_\nu(r) - J_{\nu+1}(r)$, we also have 
\begin{equation}\label{eq:J'-upper-bound}
    \sup_{2 \leq r \leq 3} \abs{J_\nu'(r)} \leq e^{9/8}(\nu+1)\frac{(3/2)^\nu}{\Gamma(\nu+1)}.
\end{equation}
Combining the two estimates \eqref{eq:J-upper-bound} and \eqref{eq:J'-upper-bound}, we obtain 
\[ \sup_{2 \leq r \leq 3} \abs{w_n(r)} + \abs{w_n'(r)} \leq C_0(\nu_n+1)\frac{(3/2)^{\nu_n}}{\Gamma(\nu_n+1)} \]
for a universal constant $C_0 > 0$ independent of $n$. Meanwhile, on $[2, 3]$, we have $\eta(r/R_n) = 1$ for all sufficiently large $n$, and thus $a_n'(r) = \chi'(r)$ and $a_n''(r) = \chi''(r)$ from which we deduce, since $\chi'$ and $\chi''$ are bounded independently of $n$, it follows that 
\begin{equation}\label{eq:delta+1-estimate-[2,3]}
    \int_2^3 \abs{2a_n'w_n' + a_n''w_n}^2dr \leq C_1(\nu_n + 1)\frac{(3/2)^{\nu_n}}{\Gamma(\nu_n+1)}
\end{equation}
for another $C_1 > 0$ independent of $n$ sufficiently large.

We now seek similar estimates on $[R_n, 2R_n]$. Using the asymptotic expansion of the Bessel function $J_\nu(r)$ as $r \rightarrow +\infty$, we note 
\begin{align*}
    w_n(r) &= \sqrt{\frac{2}{\pi}}\cos\left(r - \frac{\pi\nu_n}{2} - \frac{\pi}{4}\right) + \cO\left(\frac{\nu_n^2}{r}\right) \\
    w_n'(r) &= -\sqrt{\frac{2}{\pi}}\sin\left(r - \frac{\pi\nu_n}{2} - \frac{\pi}{4}\right) + \cO\left(\frac{\nu_n^2}{r}\right),
\end{align*}
the latter by using the recurrence relation $J_\nu'(r) = \frac{1}{2}(J_{\nu-1}(r) - J_{\nu+1}(r))$. Therefore, provided $R_n \rightarrow +\infty$ sufficiently fast, there exists a constant $M > 0$ independent of $n$ such that 
\[ \sup_{R_n \leq r \leq 2R_n} \abs{w_n(r)} + \abs{w_n'(r)} \leq M. \]
As before, on $[R_n, 2R_n]$, for sufficiently large $n$ one has $R_n > 3$ and $\chi(r) = 1$, and it follows that 
\[ a_n(r) = \eta\left(\frac{r}{R_n}\right), \quad a_n'(r) = \frac{1}{R_n}\eta'\left(\frac{r}{R_n}\right), \quad a_n''(r) = \frac{1}{R_n^2}\eta''\left(\frac{r}{R_n}\right). \]
Since $\eta$ is fixed and smooth, there exists a constant $K > 0$, depending only on $\eta$, such that 
\[ \abs{a_n'(r)} \leq \frac{K}{R_n}, \quad \abs{a_n''(r)} \leq \frac{K}{R_n^2}. \]
It follows that
\begin{align*}
    \int_{R_n}^{2R_n} \abs{2a_n'w_n' + a_n''w_n}^2dr &\leq 2\int_{R_n}^{2R_n} \abs{a_n''}^2\abs{w_n}^2dr + 8\int_{R_n}^{2R_n} \abs{a_n'}^2\abs{w_n'^2}dr \\
    &\leq 2K^2M^2\int_{R_n}^{2R_n} \frac{dr}{R_n^4} + 8K^2M^2\int_{R_n}^{2R_n} \frac{dr}{R_n^2} \vphantom{2\int_{R_n}^{2R_n} \abs{a_n''}^2\abs{w_n}^2dr + 8\int_{R_n}^{2R_n} \abs{a_n'}^2\abs{w_n'^2}dr} \\
    &= \frac{2K^2M^2}{R_n^3} + \frac{8K^2M^2}{R_n}. \vphantom{2\int_{R_n}^{2R_n} \abs{a_n''}^2\abs{w_n}^2dr + 8\int_{R_n}^{2R_n} \abs{a_n'}^2\abs{w_n'^2}dr}
\end{align*}
and conclude, since $R_n \geq 1$, that 
\begin{equation}\label{eq:delta+1-un-estimate-[R_n,2R_n]}
    \int_{R_n}^{2R_n} \abs{2a_n'w_n' + a_n''w_n}^2dr \leq \frac{C_2}{R_n} 
\end{equation}
for $C_2 > 0$ independent of $n$. Combining the two estimates \eqref{eq:delta+1-estimate-[2,3]} and \eqref{eq:delta+1-un-estimate-[R_n,2R_n]} with Stirling's approximation $\Gamma(\nu + 1) \geq \sqrt{2\pi\nu}(\nu/e)^\nu$, we find for sufficiently large $n$ and $R_n \rightarrow +\infty$ of order $R_n \gg \nu_n^2$, the estimate for \eqref{eq:delta+1-un}
\begin{equation}\label{eq:delta+1-estimate}
    \int_{B^c} \abs{(\Delta + 1)u_n}^2dx \leq \frac{C_1}{\nu_n^4} + \frac{C_2}{R_n}.
\end{equation}
We now estimate the denominator of \eqref{eq:zero-limit}. Again, expressing in polar coordinates and ignoring the surface volume of the sphere, we get 
\begin{equation}\label{eq:weighted-term-un-estimate}
    \int_{B^c} \frac{\abs{u_n}^2}{(r^2-1)^2}dx \propto \int_1^\infty \frac{\abs{a_n(r)}^2\abs{w_n(r)}^2}{(r^2-1)^2}dr.
\end{equation}
By the Weber-Schafheitlin integrals 
\begin{align*}
    \int_0^\infty \frac{J_\nu(r)^2}{r^3}dr &= \frac{1}{4\nu(\nu^2-1)}, \hspace{-3.5cm} &&\nu > 1, \\ 
    \int_0^\infty \frac{J_\nu(r)^2}{r}dr &= \frac{1}{2\nu}, \hspace{-3.5cm} &&\nu > 0,
\end{align*}
we have 
\begin{equation}\label{eq:weighted-term-estimate-0}
\begin{aligned}
    \int_1^\infty \frac{\abs{a_n(r)}^2\abs{w_n(r)}^2}{(r^2-1)^2}dr &\geq \int_3^{R_n} \frac{J_{\nu_n}(r)^2}{r^3}dr \\ 
    &= \frac{1}{4\nu_n(\nu_n^2-1)} - \int_0^3 \frac{J_{\nu_n}(r)^2}{r^3}dr - \int_{R_n}^\infty \frac{J_{\nu_n}(r)^2}{r^3}dr. 
\end{aligned}
\end{equation}
For the first remainder term, using \eqref{eq:J-upper-bound} which remains valid on $r \in (0, 3)$ gives 
\begin{equation}\label{eq:weighted-term-estimate-1}
    \int_0^3 \frac{J_{\nu_n}(r)^2}{r^3}dr \leq \frac{e^{9/4}}{2^{2\nu_n}\Gamma(\nu_n+1)^2}\int_0^3 r^{2\nu_n-3}dr = \frac{e^{9/4}3^{2\nu_n-2}}{2^{2\nu_n}(2\nu_n-2)\Gamma(\nu_n+1)^2} 
\end{equation}
while for the second remainder term 
\begin{equation}\label{eq:weighted-term-estimate-2}
    \int_{R_n}^\infty \frac{J_{\nu_n}(r)^2}{r^3}dr \leq \frac{1}{R_n^2}\int_0^\infty \frac{J_{\nu_n}(r)^2}{r}dr = \frac{1}{2\nu_nR_n^2}.
\end{equation}
Combining the three estimates \eqref{eq:weighted-term-estimate-0}, \eqref{eq:weighted-term-estimate-1}, \eqref{eq:weighted-term-estimate-2}, we conclude that 
\[ \int_1^\infty \frac{\abs{a_n(r)}^2\abs{w_n(r)}^2}{(r^2-1)^2}dr \geq \frac{1}{4\nu_n(\nu_n^2 - 1)} - \frac{e^{9/4}3^{2\nu_n-2}}{2^{2\nu_n}(2\nu_n-2)\Gamma(\nu_n+1)^2} - \frac{1}{2\nu_nR_n^2} \]
so that, again by Stirling's approximation and with $R_n \gg \nu_n$, we find for sufficiently large $n$ the estimate for \eqref{eq:weighted-term-un-estimate}
\begin{equation}\label{eq:weighted-term-estimate}
    \int_{B^c} \frac{\abs{u_n}^2}{(r^2 - 1)^2}dx \geq \frac{C_3}{\nu_n^3}
\end{equation}
for $C_3 > 0$ independent of $n$. Taking the quotient of \eqref{eq:delta+1-estimate} and \eqref{eq:weighted-term-estimate} with $R_n = \nu_n^4$ we obtain that the quotient in \eqref{eq:zero-limit} is $\cO(1/\nu_n) \rightarrow 0$ as $n \rightarrow \infty$. 
\end{proof}
The counterexample in the above argument consists of functions involving spherical harmonics of arbitrary high degree. This is the primary reason for the underlying assumption on the dimension $d \geq 2$. It would be interesting to see if the obstruction disappears, and whether one recovers the Rellich inequality on the whole space $\mathbb{R}^d \setminus \mathbb{S}^{d-1}$, when restricted to the class of functions involving spherical harmonics of degree at most a fixed $k \in \mathbb{N}_0$. It turns out to be the case for $k = 0$, which is the class of radial functions. 

\begin{proposition}\label{prop:radial-functions}
Let $d \geq 1$. Then there exist a constant $\overline{C}_d > 0$, depending only on the dimension $d$, such that
\begin{equation}\label{eq:radial-functions}
    \int_{B^c} \abs{(\Delta + 1)u}^2dx \geq \overline{C}_d \int_{B^c} \frac{\abs{u}^2}{(|x|^2-1)^2}dx. 
\end{equation}
holds for every radial function $u\in C_c^\infty(B^c)$.
\end{proposition}

\begin{proof}
Since $u = u(r)$ is radial, the Laplacian is given by
\[ \Delta u = u''(r) + \frac{d-1}{r} u'(r) \]
and making the same change of variable as before $u = r^{-\frac{d-1}{2}}v$, a direct computation (or taking $n = 0$ in \eqref{eq:delta+1-expression}) yields
\[ (\Delta + 1)u = u'' + \frac{d-1}{r} u' + u = r^{-\frac{d-1}{2}}\left(v'' + \left(1 -  \frac{(d-1)(d-3)}{4r^2}\right)v\right). \]
Set $q_d(r) = 1 - \frac{(d-1)(d-3)}{4r^2}$ and $f = v'' + q_dv$. Then in polar coordinates, it suffices to show there exists $\overline{C}_d > 0$ such that 
\begin{equation}\label{eq:radial-inequality} 
    \int_1^\infty \abs{f(r)}^2dr \geq \overline{C}_d \int_1^\infty \frac{\abs{v(r)}^2}{(r^2-1)^2}dr 
\end{equation}
where $v \in C_c^\infty((1, \infty))$ solves the initial value problem 
\[ \begin{cases} v'' + q_dv = f, &r > 1, \\
v(1) = v'(1) = 0. \end{cases} \]
% We now turn our attention to the homogeneous equation $y'' + q_dy = 0$. Let $\nu = \frac{d-2}{2}$; then since $\nu^2 - \frac{1}{4} = \frac{(d-1)(d-3)}{4}$, the equation becomes 
% \[ y'' + \left(1 - \frac{\nu^2 - 1/4}{r^2}\right)y = 0 \] 
% and it admits two (linearly independent) fundamental solutions 
% \[ y_1(r) = \sqrt{r}J_\nu(r), \quad y_2(r) = \sqrt{r}Y_\nu(r) \]
% where $Y_\nu$ denotes the second Bessel function of order $\nu$. 
As before, the homogeneous equation $y'' + q_dy = 0$ is simply the equation in the proof of Theorem \ref{thm:ext} at $n = 0$ and it admits two fundamental solutions
\[ y_1(r) = \sqrt{r}J_\nu(r), \quad y_2(r) = \sqrt{r}Y_\nu(r) \]
where $\nu = \nu_0 = \frac{d-2}{2}$ in which case $\nu^2 - \frac{1}{4} = \frac{(d-1)(d-3)}{4}$. 

Introduce their Wronskian
\[ W(r) = y_1(r)y_2'(r) - y_1'(r)y_2(r). \]
By Abel's identity, for a differential equation of the form $y'' + p(r)y' + q(r)y = 0$, the Wronskian of two solutions satisfies $W'(r) = -p(r)W(r)$ and since in the present case $p(r) = 0$, we deduce $W'(r) = 0$ and hence the Wronskian is constant. Using the asymptotic expansion of Bessel functions as $r \rightarrow \infty$ we have $\lim_{r \to \infty} W(r) = 2/\pi$ and since $W$ is constant we conclude $W(r) \equiv 2/\pi$ for all $r \geq 1$. (In particular, $y_1$ and $y_2$ are linearly independent.) 

In order to establish \eqref{eq:radial-functions}, we seek a solution of $v'' + q_dv = f$ of the form $v(r) = A(r)y_1(r) + B(r)y_2(r)$ where $A$ and $B$ satisfy 
\[ \begin{cases} A'y_1 + B'y_2 &= 0, \\ A'y_1' + B'y_2' &= f. \end{cases} \]
The determinant of this system is the Wronskian, and therefore $A(r) = -\frac{\pi}{2}y_2f$ and $B'(r) = \frac{\pi}{2}y_1f$. The initial conditions, together with $W \neq 0$, imply $A(1) = B(1) = 0$ and integrating we find 
\begin{align*}
    A'(r) &= -\frac{\pi}{2}\int_1^r y_2(s)f(s)ds, \\
    B'(r) &= \frac{\pi}{2}\int_1^r y_1(s)f(s)ds. 
\end{align*}
It follows that 
\begin{equation}\label{eq:def-G_d(r,s)}
    v(r) = \int_1^r G_d(r, s)f(s)ds, \quad G_d(r, s) = \frac{\pi}{2}(y_2(r)y_1(s) - y_1(r)y_2(s)), \quad 1 \leq s \leq r. 
\end{equation}

We claim there exists a constant $M_d > 0$, depending only on $d$ such that 
\begin{equation}\label{eq:G_d-estimate}
    \abs{G_d(r, s)} \leq M_d\min(r-s, 1), \quad 1 \leq s \leq r. 
\end{equation}
Indeed, by the asymptotic expansion of Bessel functions for fixed $\nu$ as $r \rightarrow \infty$, 
\[ J_\nu(r) = \cO(r^{-1/2}), \quad Y_\nu(r) = \cO(r^{-1/2}), \]
there exists $A_d > 0$ such that 
\[ \abs{y_j(r)} + \abs{y_j'(r)} \leq A_d \]
for all $r \geq 1$ and each $j = 1, 2$. Thus we obtain 
\[ \abs{G_d(r, s)} \leq \frac{\pi}{2}\left(\abs{y_2(r)y_1(s)} + \abs{y_1(r)y_2(s)}\right) \leq \pi A_d^2. \]
Similarly, $\abs{\partial_rG_d(r, s)} = \frac{\pi}{2} \abs{y_2'(r)y_1(s) - y_1'(r)y_2(s)} \leq \pi A_d^2$ so that, since $G_d(s, s) = 0$, the mean value theorem gives 
\[ \abs{G_d(r, s)} \leq \pi A_d^2 \abs{r - s} \]
and we deduce \eqref{eq:G_d-estimate} holds with $M_d = \pi A_d^2$.

Finally, by the representation formula for $v$ in terms of $G_d$ given in \eqref{eq:def-G_d(r,s)} and Cauchy-Schwarz, 
\begin{align*}
    \abs{v(r)}^2 &\leq \left(\int_1^r \abs{G_d(r, s)}^2ds\right)\left(\int_1^r \abs{f(s)}^2ds\right) \vphantom{\min\left(\frac{(r-1)^3}{3}, r-1\right)\norm{f}_{L^2((1, \infty))}^2} \\
    &\leq M_d^2 \left(\int_1^r \min((r-s)^2, 1)ds\right)\norm{f}_{L^2((1, \infty))}^2 \vphantom{\min\left(\frac{(r-1)^3}{3}, r-1\right)\norm{f}_{L^2((1, \infty))}^2} \\
    &\leq M_d^2 \min\left(\frac{(r-1)^3}{3}, r-1\right)\norm{f}_{L^2((1, \infty))}^2. 
\end{align*}
Dividing by $(r^2-1)^2$ and integrating over $(1, \infty)$, we return to \eqref{eq:radial-inequality}
\[ \int_1^\infty \abs{f(r)}^2dr \geq \frac{1}{KM_d^2} \int_1^\infty \frac{\abs{v(r)}^2}{(r^2-1)^2}dr, \quad K = \int_1^\infty \frac{\min((r-1)^3/3, r-1)}{(r^2-1)^2}dr. \]
Since it is readily checked the integrand is of order $r-1$ and $1/r^3$ around $r = 1$ and $r = +\infty$ respectively, we conclude $K < \infty$ which completes the proof with $\overline{C}_d = \frac{1}{KM_d^2} > 0$.
\end{proof}

\begin{remark}
In the case $d \in \{1, 3\}$ in which case $\nu^2 - \frac{1}{4} = 0$, the expression for $G_d(r, s) = \sin(r - s)$ is explicit. The inequality $\abs{\sin(t)} \leq \min(\abs{t}, 1)$ implies the proof goes through in exactly the same fashion as before and furnishes an explicit constant $\overline{C}_1 = \overline{C}_3 > 0$.
\end{remark}

Before we finish, we note again in dimensions $d \in \{1, 3\}$ where $q_d(r) = 1$ and $\Delta + 1$ acting on radial functions takes the simple form $(\Delta + 1)u = r^{-\frac{d-1}{2}}(v'' + v)$ for $u = r^{-\frac{d-1}{2}}v$ that the corresponding inequality 
\[ \int_1^\infty \abs{v''(r) + v(r)}^2dr \geq \overline{C}_d \int_1^\infty \frac{\abs{v(r)}^2}{(r^2-1)^2}dr \]
can be established by a simple and direct proof which follows from the classical weighted Hardy inequality \cite{MR4089473}
\begin{equation}\label{eq:classical-weighted-hardy}
    \int_\bbR \abs{f'}^2r^\alpha \geq \frac{(\alpha-1)^2}{4}\int_\bbR \abs{f}^2r^{\alpha-2}dr 
\end{equation}
valid for all $\alpha \in \bbR$ and $f \in C_c^\infty(\bbR \setminus \{0\})$. Indeed, we may factorise the operator $\Delta + 1$ into the following product of first order operators 
\[ \frac{d^2}{dr^2} + 1 = \left(\frac{d}{dr} + i\right)\left(\frac{d}{dr} - i\right). \]
Let $g \vcentcolon= (\frac{d}{dr} - i)f$ and $\tilde{g} = e^{ir}g$, and let $\tilde{f} = e^{-ir}f$. Then by direct differentiation we have 
\[ \left(\frac{d^2}{dr^2}+1\right)f = \left(\frac{d}{dr} + i\right)g = e^{-ir}\tilde{g}' \,\,\text{and} \left(\frac{d}{dr} - i\right)f = e^{ir}\tilde{f}'. \]
Applying \eqref{eq:classical-weighted-hardy} with $\alpha = 0$ and $\alpha = -2$ sequentially gives 
\begin{align*}
    \int_\bbR \abs{f'' + f}^2dr = \int_\bbR \smallabs{e^{-ir}\tilde{g}'}^2dr \geq \frac{1}{4} \int_\bbR \abs{g}^2r^{-2}dr &= \frac{1}{4} \int_\bbR \abs{\left(\frac{d}{dr} - i\right)f}^2r^{-2}dr \\
    &= \frac{1}{4}\int_\bbR \smallabs{e^{ir}\tilde{f}'}^2 r^{-2} dr \vphantom{\frac{1}{4} \int_\bbR \abs{\left(\frac{d}{dr} - i\right)f}^2r^{-2}dr} \\
    &\geq \frac{9}{16} \int_\bbR \smallabs{\tilde{f}}^2r^{-4}dr \vphantom{\frac{1}{4} \int_\bbR \abs{\left(\frac{d}{dr} - i\right)f}^2r^{-2}dr}\\ 
    &= \frac{9}{16} \int_\bbR \abs{f}^2r^{-4}dr. \vphantom{\frac{1}{4} \int_\bbR \abs{\left(\frac{d}{dr} - i\right)f}^2r^{-2}dr}
\end{align*}
Translating $r \mapsto r - a$, for each fixed $a \in \bbR$ we obtain
\[ \int_\bbR \abs{f'' + f}^2dr \geq \frac{9}{16}\int_\bbR \abs{f}^2(r-a)^{-4}dr \]
for $f \in C_c^\infty(\bbR \setminus \{a\})$. 
Adding the inequalities for $a = 1$ and $a = -1$, along with the arithmetic-mean geometric-mean inequality, yields the inequality \eqref{eq:radial-inequality} for $f \in C_c^\infty(\bbR \setminus \{\pm 1\})$ with $\overline{C}_d = \frac{9}{16}$ after taking $f = v$ and extending by zero outside of $(1, \infty)$.

\subsection*{Acknowledgements}
We are grateful to Ari Laptev for sharing the problem and for subsequent discussions. The third author is grateful to the Hausdorff Research Institute for Mathematics for their support and hospitality where part of the work on this paper was undertaken. 

\printbibliography

\end{document}